\documentclass{amsart}

\usepackage{amsmath, amssymb, amsthm} 

\theoremstyle{plain}                       
\newtheorem{theorem}{Theorem}          
\newtheorem{lemma}[theorem]{Lemma}         
\newtheorem{corollary}[theorem]{Corollary} 

\theoremstyle{remark}                      
\newtheorem*{acknow}{Acknowledgements}          

\newcommand\normcl[1]{\overline{#1}^{\|\cdot\|}}

\newcommand\pair[1]{\langle\, #1\,\rangle}

\DeclareMathOperator\supp{supp}

\DeclareMathOperator\dom{dom}
\DeclareMathOperator\linsp{span}

\newcommand{\lt}{L}                        
\newcommand{\rt}{R}                        

\newcommand\sm{\setminus}
\newcommand\col{\colon}
\newcommand\sub{\subseteq}

\newcommand\cstarred{\mathrm{C}^*_r}
\newcommand\cstar{\mathrm{C}^*}

\newcommand{\R}{\mathbf{R}}                
\newcommand{\complex}{\mathbf{C}}          

\newcommand{\A}{\mathrm{A}}              
\newcommand{\B}{\mathrm{B}}              
\newcommand{\posdef}{\mathrm{P}}

\newcommand{\set}[2]{\{\,{\textstyle#1};\,{\textstyle #2}\,\}}

\newcommand{\inv}{^{-1}}                   
\newcommand\conj{\overline}

\newcommand{\st}{\mathrm{st}} 

\newcommand{\cop}{\Gamma}
\newcommand{\id}{\mathrm{id}}

\newcommand\G{\mathbb{G}}
\renewcommand\H{\mathbb{H}}      
\newcommand\K{\mathbb{K}}

\newcommand\lone{\mathrm{L}^1}
\newcommand\ltwo{\mathrm{L}^2}
\newcommand\linfty{\mathrm{L}^\infty}
\newcommand{\vn}{\mathrm{VN}}
\newcommand{\C}{\mathrm{C}}

\newcommand{\M}{\mathrm{M}}

\newcommand\vntimes{\mathop{\overline\otimes}}

\newcommand\dual[1]{\widehat{#1}} 
\newcommand\latedual{\widehat{\phantom{x}}} 

\newcommand\sw[1]{_{#1}}

\begin{document}

\title[Quantum subgroups, invariant subalgebras]%
      {Compact quantum subgroups and\\ 
       left invariant C*-subalgebras of\\
       locally compact quantum groups}

\author{Pekka Salmi}

\address{University of Oulu, Department of Mathematical Sciences,
PL 3000, FI-90014 Oulun yliopisto, Finland}

\email{pekka.salmi@iki.fi}

\thanks{Research supported by Academy of Finland}

\begin{abstract}
We show that there is a one-to-one correspondence between compact quantum 
subgroups of a co-amenable locally compact quantum group $\G$ and certain 
left invariant C*-sub\-alge\-bras of $\C_0(\G)$. We also prove that
every compact quantum subgroup of a co-amenable quantum group is co-amenable. 
Moreover, there is a one-to-one correspondence between open subgroups 
of an amenable locally compact group $G$ and non-zero, invariant 
C*-subalgebras of the group C*-algebra $\cstar(G)$.
\end{abstract}


\maketitle

\section{Introduction} 

As shown by Lau and Losert \cite{lau-losert:complemented},
there is a one-to-one correspondence between 
compact subgroups $H$ of a locally compact group $G$ and
non-zero, left translation invariant C*-subalgebras $X$ of $\C_0(G)$,
the continuous functions vanishing at infinity. The correspondence is given by 
\begin{gather*}
X = \C_0(G/H) = \set{f\in \C_0(G)}{\rt_s f = f\text{ for every } s\in H}\\
H = \set{s\in G}{\rt_s f = f\text{ for every } f\in X}
\end{gather*} 
where $\rt_s$ denotes the right translation operator on $\C_0(G)$ 
by $s$ in $G$. The normality of the subgroup is also characterised in terms 
of the subalgebra: $H$ is normal if and only if $X$ is right translation 
invariant. Preceding Lau and Losert, de Leeuw and Mirkil 
\cite{de-leeuw-mirkil:invariant-subalgebras} discovered the correspondence 
of compact subgroups and translation invariant C*-subalgebras 
in the case of locally compact abelian groups. The result of Lau and 
Losert is also related to a group of results due to Takesaki and 
\hbox{Tatsuuma}~\cite{takesaki-tatsuuma:duality-subgroups}
concerning invariant subalgebras of $\linfty(G)$, 
the group von Neumann algebra $\vn(G)$, the group algebra $\lone(G)$
and the Fourier algebra $\A(G)$. 

The purpose of the present paper is to consider the result of Lau and Losert 
in the context of locally compact quantum groups as defined by 
Kustermans and Vaes~\cite{kustermans-vaes:lcqg}. We show that compact 
quantum subgroups of a co-amenable locally compact quantum group $\G$ 
correspond to certain left invariant C*-subalgebras of 
$\C_0(\G)$. The construction of the compact quantum subgroup from a 
left invariant C*-subalgebra is quite general, 
but the left invariant C*-subalgebras arising from 
compact quantum subgroups always have certain properties
that we need to assume to obtain uniqueness. 

The notion of left invariant C*-subalgebra used here is analogous 
to the notion of left co-ideal.  
In \cite{tomatsu:coideal}, Tomatsu studied the relation between
compact quantum subgroups and right co-ideals for co-amenable compact
quantum groups, using the von Neumann algebraic setting. 
Although there is some overlap between the present paper
and \cite{tomatsu:coideal}, the techniques are quite different 
and the context of \cite{tomatsu:coideal} more restricted.
Co-ideals and compact quantum subgroups are also related to 
idempotent states. Franz and Skalski \cite{franz-skalski:idempotent-states} 
studied these relations recently for co-amenable compact quantum groups.
Finally, already Enock \cite{enock:sous-facteur}
gave a version of the Takesaki--Tatsuuma duality for 
von Neumann algebraic quantum groups.

As a necessary intermediate result, we  prove that a compact quantum 
subgroup of a co-amenable quantum group is also co-amenable. 
This result encompasses the fact that the quotient group
of an amenable locally compact group by a normal, open subgroup is amenable. 
Tomatsu \cite{tomatsu:coideal} proved this result 
in the special case when also the ambient quantum group is compact.

Finally, we  consider the particular case 
of co-amenable, co-commutative quantum groups, 
given by group C*-algebras $\cstar(G)$ of amenable locally compact groups.
In this setting we  show the dual version of the Lau and Losert 
result: there is a one-to-one correspondence between 
open subgroups of $G$ and non-zero, invariant C*-subalgebras of $\cstar(G)$.

\section{Preliminaries}
\label{sec:prel}

We define the notion of locally compact quantum group 
implicitly through its reduced C*-algebra. 
The reduced C*-algebraic approach to quantum groups 
is due to Kustermans and Vaes \cite{kustermans-vaes:lcqg}.
A (\emph{locally compact}) \emph{quantum group} $\G$
is determined by the following structures:
\begin{itemize}
\item a C*-algebra $\C_0(\G)$ (the \emph{reduced C*-algebra} 
      associated with $\G$)
\item a \emph{co-multiplication} $\cop$ on $\C_0(\G)$
\item \emph{left and right Haar weights} $\phi$ and $\psi$ on $\C_0(\G)$.
\end{itemize}
The co-multiplication $\cop$ is a non-degenerate $*$-homo\-mor\-phism 
from $\C_0(\G)$ to the multiplier algebra $\M(\C_0(\G)\otimes \C_0(\G))$ 
of the minimal C*-algebraic tensor product $\C_0(\G)\otimes \C_0(\G)$ such that
\[
(\id\otimes \cop)\cop = (\cop\otimes \id)\cop \qquad\text{(co-associativity)}.
\]
We suppose further that the linear spans of
$\cop(\C_0(\G))(\C_0(\G)\otimes 1)$ 
and $\cop(\C_0(\G))(1\otimes \C_0(\G))$ are dense in $\C_0(\G)\otimes
\C_0(\G)$. The left Haar weight $\phi$ is a faithful  
KMS-weight on $\C_0(\G)$ that is \emph{left invariant}
in the sense that
\[
\phi\bigl((\mu\otimes \id)\cop(a)\bigr) = \mu(1)\phi(a)
\]
for every positive functional $\mu$ on $\C_0(\G)$ and 
for every positive $a$ in $\C_0(\G)$ such that $\phi(a)<\infty$.
The right Haar weight is defined analogously,
using right slices $\id\otimes \mu$ instead of left.
For more details, especially in regard to 
the weights, see \cite{kustermans-vaes:lcqg}.
We shall denote the multiplier algebra of $\C_0(\G)$ by $\C(\G)$.

In the commutative case, that is, when $\C_0(\G)$ is commutative, 
there is an actual locally compact group $G$ such that 
$\C_0(\G) = \C_0(G)$ -- the C*-algebra of 
continuous functions on $G$ vanishing at infinity. 
In this case the co-multiplication is defined by
\[
\cop(f)(s,t) = f(st) \qquad(f\in \C_0(G),\, s,t\in G).
\]
The left and the right Haar measures of $G$ give the 
left and the right Haar weights via integration. 

We should specify one detail about the definition of quantum group
because it is needed elsewhere in the paper.
The maps appearing in the identities 
defining co-asso\-cia\-tiv\-ity and left invariance are, 
in fact, strict extensions to the appropriate multiplier algebras. 
The strict topology on the multiplier algebra $\M(A)$ 
of a C*-algebra $A$ is defined by the seminorms 
$u\mapsto \|ua\| + \|au\|$ where $a$ runs through the elements of $A$. 
Certain maps admit unique extensions from a C*-algebra to its 
multiplier algebra such that the extension is strictly continuous.
(In the literature these extensions are often noted to 
be strictly continuous on \emph{bounded sets} 
but that restriction is in fact unnecessary. 
As shown by Taylor~\cite[Corollary~2.7]{taylor:strict},
the strict topology coincides with the so-called
bounded strict topology -- the 
strongest locally convex topology that agrees 
with the strict topology on norm-bounded sets.
It follows that if a linear map from $\M(A)$ to a locally convex space
is strictly continuous on bounded sets, it is strictly continuous.)
For example, all bounded functionals have such \emph{strict extensions}.
More generally, consider C*-algebras $A$ and $B$, 
and let $\mu\in A^*$.
Then also the slice map $\mu\otimes \id\col A\otimes B\to B$ 
has a strict extention to a map from $\M(A\otimes B)$ to $\M(B)$
(see for example \cite[Section~7]{kustermans:one-parameter} or 
\cite[Appendix~A]{mas-nak-wor}).
Non-degenerate completely positive maps (in particular, 
non-degenerate $*$-homomorphisms) form  another important 
class of functions admitting strict extensions
\cite[Corollary~5.7]{lance:hilbert-cstar}. 
A completely positive map $\theta\col A \to \M(B)$ is said 
to be \emph{non-degenerate} if, for some bounded approximate 
identity $(e_i)$ in $A$, the net $(\theta(e_i))$
converges strictly to the identity in $\M(B)$.
We shall often use strict extensions without further mention.

Let $\G$ be a quantum group. The Gelfand--Naimark--Segal 
construction applied to the left Haar weight $\phi$ of $\G$
gives a Hilbert space $\ltwo(\G)$ and a faithful representation 
of $\C_0(\G)$ on $\ltwo(\G)$. We shall identify $\C_0(\G)$ with 
its image under this representation and consider $\C_0(\G)$ as 
a C*-subalgebra of $\B(\ltwo(\G))$, the bounded operators on $\ltwo(\G)$. 
There is a unitary operator $W$ in $\B(\ltwo(\G)\otimes \ltwo(\G))$ 
that determines the co-multiplication of $\G$ by
\[
\cop(a) = W^*(1\otimes a)W \qquad(a\in \C_0(\G))
\]
and satisfies the pentagonal equation
\[
W_{12}W_{13}W_{23} = W_{23}W_{12}.
\]
Here $W_{12} = W\otimes 1$, $W_{23} = 1\otimes W$ and
$W_{13} = (\Sigma\otimes \id)(W\otimes 1)(\Sigma\otimes \id)$ where
$\Sigma$ denotes the flip map on $\ltwo(\G)\otimes \ltwo(\G)$. 
The operator $W$ is called the \emph{left multiplicative unitary}. 
It is in fact an element of 
$\M(\C_0(\G)\otimes \B_0(\ltwo(\G)))$, where $\B_0(\ltwo(\G))$ 
denotes the C*-algebra of compact operators on $\ltwo(\G)$.  

Let $\B(\ltwo(\G))_*$ denote the weak*-continuous functionals on 
$\B(\ltwo(\G))$. Then 
\[
\C_0(\G) = \normcl{\set{(\id\otimes\omega)W}{\omega\in \B(\ltwo(\G))_*}}.
\]
On the other hand, also 
\[
\C_0(\dual \G) = \normcl{\set{(\omega\otimes\id)W}{\omega\in \B(\ltwo(\G))_*}}
\]
is a C*-subalgebra of $\B(\ltwo(\G))$. 
Moreover, $W\in \M(\C_0(\G)\otimes \C_0(\dual\G))$. 
Put $\dual W = \Sigma W^* \Sigma$ and  define 
\[
\dual{\cop}(x) = \dual W^* (1\otimes x)\dual W \qquad(x\in\C_0(\dual \G)).
\]
One of the main results of \cite{kustermans-vaes:lcqg} is that 
$\dual\G$ is a locally compact quantum group --
the \emph{dual quantum group} of $\G$ -- and 
the analogue of the Pontryagin duality theorem holds: 
$\dual{\dual \G} = \G$.

As an example, consider the commutative case when $\G = G$ 
for some locally compact group $G$. The operators $W$ and $\dual W$ 
on $\ltwo(G\times G)\cong \ltwo(G)\otimes \ltwo(G)$ are given by 
\[
\begin{split}
W\xi(s,t) &= \xi(s,s\inv t)\\
\dual W\xi(s,t) &= \xi(ts,t)
\end{split}
\qquad(\xi\in\ltwo(G\times G),\, s,t\in G)
\]
(the identities are understood to hold almost everywhere). 
Let $\lambda$ be the left 
regular representation of $G$ so that $\lambda(f)$ is the convolution
operator on $\ltwo(G)$ determined by $f$ in $\lone(G)$. 
Moreover, let $M_g$ denote the operator of pointwise multiplication 
by $g$ in $\C_0(G)$. 
For fixed $\xi$ and $\zeta$ in $\ltwo(G)$, define $\omega$ in $\B(\ltwo(G))_*$
by $\omega(x) = \pair{x\xi \mid \zeta}$,  $x\in \B(\ltwo(G))$.
Then $(\id\otimes \omega)W = M_{\conj{\zeta}*\check\xi}$
(where $\check\xi(s)=\xi(s\inv)$) 
and $(\omega\otimes \id) W = \lambda(\xi\conj\zeta)$.
We see that $\C_0(\dual\G) = \cstarred(G)$ --
the reduced group C*-algebra of $G$.
The co-multiplication of $\cstarred(G)$ is defined by
\[
\cop(\lambda(s)) = \lambda(s)\otimes \lambda(s) \qquad(s\in G),
\]
and the Plancherel weight \cite[Section VII.3]{takesaki:vol2}
acts as both the left and the right Haar weight.
This case is called the \emph{co-commutative case} because
these are exactly the quantum groups $\G$ 
such that the co-multiplication 
satisfies the co-commutativity condition $\sigma \cop = \cop$,
where $\sigma$ is the flip map defined
by $\sigma(a\otimes b) = b\otimes a$ ($a,b\in \C_0(\G)$).

Let $\G$ again be an arbitrary quantum group. 
We shall also use the von Neumann algebraic side of $\G$.
Let $\linfty(\G)$ denote the von Neumann algebra generated 
by $\C_0(\G)$ in $\B(\ltwo(\G))$.  
Then the co-multiplication extends to a normal $*$-homormorphism
$\cop\col \linfty(\G)\to \linfty(\G)\vntimes\linfty(\G)$
that is still co-associative. Also the Haar weights have natural 
counterparts on the von Neumann algebra. 
The reader is referred to \cite{kustermans-vaes:lcqg-vn, vandaele:lcqg}
for a detailed account on the von Neumann algebraic theory of 
quantum groups. 
Let $\lone(\G)$ be the closed subalgebra of $\C_0(\G)^*$ consisting 
of the functionals in $\B(\ltwo(\G))_*$ restricted to $\C_0(\G)$,
so $\lone(\G)$ can be identified with the predual of $\linfty(\G)$.

Still another facet of $\G$ is given by 
the universal C*-algebraic approach due to 
Kustermans \cite{kustermans:universal}.
At this point, we need the antipode operator, which plays the role 
of inverse. The antipode $S$ associated with the quantum 
group $\G$ is a closed, densely defined operator on $\C_0(\G)$.
The elements $(\id\otimes \sigma)W$, $\sigma\in\B(\ltwo(\G))_*$,
form a core of $S$ and 
$S\bigl((\id\otimes\sigma)W\bigr) = (\id\otimes\sigma)W^*$
\cite[Proposition~8.3]{kustermans-vaes:lcqg}.
We denote the domain of $S$ by $\dom S$.
For every $\mu$ in $\C_0(\G)^*$, define $\conj\mu$ in $\C_0(\G)^*$ 
by $\conj\mu(a) = \conj{\mu(a^*)}$,  $a\in \C_0(\G)$.
Following \cite{kustermans:universal}, we define
\[
\lone_*(\G) = \set{\omega\in\lone(\G)}{\exists \eta\in\lone(\G)
  \text{ such that }\conj{\omega}(S(x)) = \eta(x) \;\forall x\in \dom(S)}.
\]
For every $\omega$ in $\lone_*(\G)$, write $\omega^*$ for the
functional $\eta$ associated with $\omega$ in the definition of $\lone_*(\G)$. 
The map $\omega\mapsto \omega^*$ is an involution on $\lone_*(\G)$, 
and $\lone_*(\G)$ is a Banach $*$-algebra with respect to the 
norm $\|\omega\|_* = \max\{\|\omega\|, \|\omega^*\|\}$. 
The universal C*-algebra $\C^u_0(\G)$ is
the universal C*-completion of the 
Banach $*$-algebra $\lone_*(\dual\G)$ \cite[p. 296]{kustermans:universal}.
Also $\C^u_0(\G)$ has a co-multiplication, which 
we denote by $\cop^u$.
There is always a surjective $*$-homomorphism 
$\rho\col\C^u_0(\G)\to\C_0(\G)$ such that 
$(\rho\otimes\rho)\cop^u = \cop\rho$.
In the commutative case $\C^u_0(\G) = \C_0(\G)$,
but in the co-commutative case $\C^u_0(\G)$ is the 
universal group C*-algebra of $G$, which is different from 
the reduced group C*-algebra whenever $G$ is non-amenable. 
In general, $\C_0^u(\G) = \C_0(\G)$ exactly when 
$\G$ is \emph{co-amenable}, that is,
when there is $\epsilon$ in $\C_0(\G)^*$ such that 
\[
(\epsilon\otimes \id)\cop = (\id\otimes \epsilon)\cop = \id.
\]
The functional $\epsilon$ is called the \emph{co-unit};
it is a $*$-homomorphism. 
In this paper we shall concentrate on co-amenable quantum groups,
and $\epsilon$ denotes always the co-unit.
See \cite{bedos-tuset:amenability} for several characterisations of
co-amenability.

A quantum group $\G$ is said to be \emph{compact} if $\C_0(\G)$ is unital
(so $\C(\G) = \C_0(\G)$).
In this case, the left Haar weight is not just a weight but a positive 
functional. Normalisation of the left Haar weight gives a state, which
is both left and right invariant. This state is called the 
\emph{Haar state} of~$\G$. The theory of compact quantum groups, 
developed by Woronowicz 
\cite{woronowicz:compact-matrix-pseudogroups,woronowicz:compact},
predates the theory of locally compact quantum groups. 
It is also more satisfying because the existence of 
Haar weights need not be assumed, but as Woronowicz showed, the Haar state
can be constructed. (Faithfulness of the constructed state is not, however,
guaranteed.) There is a very nice survey on compact quantum groups by 
Maes and Van~Daele~\cite{maes-vandaele:compact}.

Let $\G$ be a quantum group. We say that $(\H, \pi)$ 
is a \emph{closed quantum subgroup} of $\G$ 
if $\H$ is a quantum group and $\pi\col \C_0^u(\G)\to \C_0^u(\H)$ 
is a surjective $*$-homo\-mor\-phism 
such that $(\pi\otimes \pi)\cop^u = \cop_{\H}^u \pi$. 
The motivation for the definition is of course the case of a 
closed subgroup $H$ of a locally compact group $G$; then 
$\pi$ is the restriction map from $\C_0(G)$ onto $\C_0(H)$.
Recall that our framework will be that of co-amenable 
quantum groups and we shall only consider compact quantum subgroups, 
in which case the morphism $\pi$ goes from $\C_0(\G)$ to $\C(\H) = \C_0(\H)$
(see section~\ref{sec:coamen}). 
We say that two closed quantum subgroups $(\H, \pi)$ and $(\H', \pi')$ 
are \emph{isomorphic} if there is a $*$-isomorphism 
$\rho\col \C_0^u(\H)\to \C_0^u(\H')$ such that $\rho\pi = \pi'$. In this case 
we write $(\H, \pi)\cong(\H',\pi')$.
There are also other definitions for closed quantum subgroup.
In \cite{vaes-vainerman:low-dimensional} closed quantum subgroup is 
defined on the von Neumann algebraic level, and 
in \cite{vaes:imprimitivity} the definition goes through 
both the universal level and the von Neumann algebraic one.
It is not quite clear when the different notions coincide.

We also need an analogue of left translation and left invariance. 
The dual space $\C_0(\G)^*$ is a Banach algebra under the multiplication
\[
\mu * \nu = (\mu\otimes \nu)\cop\qquad(\mu,\nu\in \C_0(\G)^*).
\]
In the commutative case $\G = G$, the above multiplication
is the convolution on the measure algebra $\M(G)$. 
In the co-commutative case $\G = \dual G$, 
the above multiplication is the pointwise product on
the reduced Fourier--Stieltjes algebra $\B_r(G)$.

To imitate the notion of left translation, define
\[
\lt_\mu(a) = (\mu\otimes \id)\cop(a)\qquad(\mu\in \C_0(\G)^*, a\in \C_0(\G)).
\]
In the commutative case $\G= G$, the operator $\lt_{\delta_s}$, 
determined by the point mass $\delta_s$ at $s$ in $G$, is exactly the 
left translation operator $\lt_s$. In general, the map 
$(\mu,a)\mapsto \lt_\mu(a)$ is a right action of the Banach algebra 
$\C_0(\G)^*$ on $\C_0(\G)$. 
A subspace $X\sub \C_0(\G)$ is said to be \emph{left invariant} 
if the operators $\lt_\mu$ leave $X$ fixed, i.e., 
\[
\lt_\mu(x) = (\mu\otimes \id)\cop(x) \in X 
\qquad\text{for every $\mu$ in $\C_0(\G)^*$ and $x$ in $X$.}
\]
As noted in \cite{salmi:quantum-luc}, if $\G$ is 
co-amenable, then the C*-algebra $\C_0(\G)$ has the slice map property 
introduced by Wassermann \cite{wassermann:first-slice}. 
It follows that every left invariant C*-subalgebra $X$ of 
co-amenable $\C_0(\G)$ satisfies 
\[
\cop(x) \in \M(\C_0(\G)\otimes X)\qquad\text{for every $x$ in } X.
\]
We shall not use the above property explicitly, but it is 
analogous to the notion of left co-ideal. 

We end this section with a word on notation. 
Now let $A$ be a C*-algebra.
It follows from Cohen's factorisation theorem that, 
for every $\mu$ in $A^*$, 
there exist $\nu$ and $\eta$ in $A^*$ and
$a$ and $c$ in $A$ such that 
\[
\mu(b) = \nu(ab) =: \nu.a(b) \qquad(b\in A)
\]
and
\[
\mu(b) = \eta(ac) =: c.\eta(b) \qquad(b\in A).
\] 
The factorisation applies in particular
to weak*-continuous functionals on $\B(H)$ 
because  $\B(H)_*=\B_0(H)^*$. 

\section{Compact quantum subgroups from left invariant C*-subalgebras}
\label{sec:subalgebra->subgroup}

In this section we construct a compact quantum subgroup 
from a left invariant C*-subalgebra of a quantum group. 
Let $\G$ be a co-amenable quantum group, and let 
$X$ be a non-zero, left invariant C*-subalgebra of $\C_0(\G)$. 

We say that a non-degenerate $*$-homomorphism $\rho\col \C_0(\G)\to\M(A)$,
where $A$ is a C*-algebra, is \emph{$X$-trivial} if
\[
\rho(x) = \epsilon(x)1_{\M(A)}\qquad\text{for every $x$ in }X.
\]
Let $J_X = \bigcap_\rho \ker\rho$
where the intersection is taken over the equivalence classes of 
non-degenerate $X$-trivial representations of $\C_0(\G)$.
Then $J_X$ is a closed ideal in $\C_0(\G)$.

\begin{lemma} \label{lemma:strict-J_X}
An element $a$ in $\C(\G)$ is in the strict closure of $J_X$ 
if and only if\/ $\widetilde\rho(a) = 0$ 
for every non-degenerate, $X$-trivial $*$-homomorphism
$\rho\col \C_0(\G)\to\M(A)$
\textup{(}here $\widetilde\rho$ denotes the strict extension\textup{)}. 
\end{lemma}

\begin{proof}
Suppose that the latter condition holds for $a$ in $\C(\G)$.
Let $(e_i)$ be a bounded
approximate identity of $\C_0(\G)$ so that $ae_i\to a$ strictly.
Now if $\rho$ is a non-degenerate $X$-trivial representation
of $\C_0(\G)$, then $\rho(ae_i) = \widetilde\rho(a)\rho(e_i) = 0$.
It follows that $ae_i\in J_X$ and so $a$ is in the strict closure of $J_X$.

The converse is clear because 
every non-degenerate, $X$-trivial $*$-homomorphism
$\rho\col \C_0(\G)\to\M(A)$
vanishes on $J_X$ (just combine $\rho$ with a non-degenerate, 
faithful representation of $A$)
and the extension $\widetilde\rho$ is strictly continuous.
\end{proof}

\begin{theorem} \label{thm:subalgebra->subgroup}
Suppose that $\G$ is a co-amenable quantum group and that $X$ is a 
non-zero, left invariant C*-subalgebra of\/ $\C_0(\G)$. Then there 
is a compact quantum subgroup $(\H,\pi)$ of\/ $\G$ such that 
$\C(\H) = \C_0(\G)/J_X$ and $\pi\col \C_0(\G)\to \C_0(\G)/J_X$ 
is the quotient map.
\end{theorem}

\begin{proof}
We begin by showing that the quotient map $\pi\col \C_0(\G)\to \C_0(\G)/J_X$ 
is $X$-trivial. Let $x\in X$. 
Then the strict extension of every $X$-trivial non-degenerate
$*$-homomorphism vanishes at $x-\epsilon(x)1$, 
so $x-\epsilon(x)1$ is in the strict closure of $J_X$
by Lemma~\ref{lemma:strict-J_X}.
Therefore $x-\epsilon(x)1$ is in the kernel of the strict extension of $\pi$. 
It follows that  $\pi(x) = \epsilon(x)1$ so $\pi$ is $X$-trivial. 

Let $x\in X$ be non-zero so that $\mu(x) = 1$ for some $\mu\in\C_0(\G)^*$.
Then $y = (\mu\otimes\id)\cop(x)$ is in $X$ by left invariance,
and $\epsilon(y) = 1$.
Since $\pi$ is $X$-trivial, $\pi(y) = 1$
and hence $\C_0(\G)/J_X$ is in fact unital. 

Denote the unital C*-algebra $\C_0(\G)/J_X$ tentatively by $\C(\H)$.
We want to define a co-multiplication $\cop_\H$ on $\C(\H)$ by
\[
\cop_\H(\pi(a)) = (\pi\otimes \pi)\cop(a)\qquad(a\in \C_0(\G)).
\]
For this definition to make sense, we need that 
$(\pi\otimes \pi)\cop(a) = 0$ whenever $a\in J_X$.
In other words, we should show that 
$(\pi\otimes\pi)\cop$ is $X$-trivial. 
Let $x\in X$. Now for every $\mu$ in $\C_0(\G)$, 
\[
\pi\bigl((\mu\otimes \id)\cop(x)\bigr) = \mu(x)1
\]
because $X$ is left invariant and $\pi$ is $X$-trivial.
Therefore
\[
(\pi\otimes\pi)\cop(x) =(\pi\otimes \id)(x\otimes 1) = \epsilon(x)(1\otimes 1),
\]
so $(\pi\otimes\pi)\cop$ is $X$-trivial as required. 

It follows from the definition of $\cop_{\H}$ and 
the properties of $\cop$ that $\cop_{\H}$ is
a co-multiplication such that the linear spans 
of $\cop_{\H}(\C(\H))(\C(\H)\otimes 1)$ and $\cop_{\H}(\C(\H))(1\otimes\C(\H))$ 
are dense in $\C(\H)\otimes \C(\H)$.
Since $\C(\H)$ is unital, there is a state $\phi_{\H}$ of $\C(\H)$ that is 
both left and right invariant (see for example  
\cite{maes-vandaele:compact}). 
It is not clear, however, that $\phi_{\H}$ is faithful,
but we can resolve this issue by taking a further quotient:
as shown in \cite[Proposition 5.4.8]{timmermann:quantum-groups},
there is a compact quantum group $\K$
and a surjective $*$-homomorphism $\rho\col \C(\H)\to \C(\K)$ 
such that $(\rho\otimes\rho)\cop_{\H} = \cop_{\K}\rho$.
By Theorem~\ref{thm:co-amenable-subgroup} below,
$\K$ is co-amenable, so $\C(\K)$ is isomorphic with the universal 
C*-algebra $\C^u(\K)$ by \cite[Theorem~3.1]{bedos-tuset:amenability}.
It follows from \cite[Proposition 5.4.8]{timmermann:quantum-groups}
that $\C(\K) \cong \C(\H) \cong \C^u(\K)$, so $\H = \K$ is a compact quantum 
subgroup of~$\G$.
\end{proof}

With the exception of the very last step, 
the preceding construction works also when $\G$ 
is not necessarily co-amenable, in which case the construction
is applied to the universal C*-algebra $\C_0^u(\G)$
and the co-unit of $\C_0^u(\G)$. 
The resulting quotient C*-algebra is \emph{a} C*-algebra
associated with a compact quantum group, but it is not 
clear whether it is the \emph{universal} C*-algebra of a 
compact quantum group
(which is necessary for the definition of closed quantum subgroup 
to be satisfied).

It turns out that different left invariant C*-subalgebras
may induce the same compact quantum subgroup through the
preceding construction. 
We shall next introduce a condition that is necessary for 
a one-to-one correspondence result.
We say that a left invariant C*-subalgebra $X$ is \emph{symmetric} if
\[
W(x\otimes 1)W^* \in \M(X\otimes \B_0(H))\qquad(x\in X).
\]
Perhaps a more descriptive term would be \emph{co-action symmetric},
following \cite{tomatsu:coideal}, 
but as co-actions are not prominently present, 
let us use the simpler terminology.
There should be no risk of confusion, although 
in the context of Kac algebras 
the term ``symmetric'' is sometimes used instead of ``co-commutative''.
In the case of classical groups, every left invariant C*-subalgebra is 
automatically symmetric due to commutativity. By Theorem~\ref{thm:dual-normal} 
below, an invariant C*-subalgebra of a co-amenable, co-commutative 
quantum group is symmetric if and only if the corresponding open subgroup 
is normal. We shall see in section~\ref{sec:subgroup->subalgebra} 
that the left invariant 
C*-subalgebra arising from a compact quantum subgroup is necessarily symmetric. 
Vaes and Vainerman \cite{vaes-vainerman:low-dimensional}
used a similar condition to define \emph{normality}
for closed quantum subgroups in the von Neumann algebraic setting.

We shall use the symmetry property to characterise the dual space of 
$\C(\H)$ when $\H$ is the compact quantum subgroup of $\G$ associated
with a symmetric left invariant C*-subalgebra $X$.
Denote the state space of $\C_0(\G)$ by $\st(\C_0(\G))$ 
(that is, $\st(\C_0(\G))$ is the collection of all positive
functionals on $\C_0(\G)$  with norm~$1$). 
Define
\[
F_0 = \set{\mu\in \st(\C_0(\G))}%
{(\id\otimes \mu)\cop(x)= x \text{ for every }x\in X}.
\]
We shall proceed to show that when $X$ is symmetric, $F_0$ consists 
precisely of those states that have $X$-trivial GNS-representations. 

\begin{lemma} 
\label{lemma:basic-F-0}
\begin{enumerate}
\item\label{item:F-0} 
$F_0 = \set{\mu\in\st(\C_0(\G))}{\mu = \epsilon\text{ on $X$}}$.
\item\label{item:mult}
 If $\mu\in F_0$, then 
    $\mu(ax) = \mu(a)\mu(x)$ and $\mu(xa) = \mu(x)\mu(a)$
       for every $a$ in $\C_0(\G)$ and $x$ in $X$.
\end{enumerate}
\end{lemma} 

\begin{proof}
(\ref{item:F-0})
For every $\mu$ in $F_0$ and $x$ in $X$,
\[
\mu(x) = (\epsilon\otimes\mu)\cop(x) = \epsilon(\id\otimes\mu)\cop(x)
       = \epsilon(x).
\]
The converse follows from the left invariance of $X$ and
the identity $(\id\otimes\epsilon)\cop = \id$.

(\ref{item:mult}) 
Let $\mu\in F_0$ and $x\in X$. By (\ref{item:F-0}), 
\[
\mu(x^* x) = \epsilon(x^* x) = \epsilon(x^*)\epsilon(x) = \mu(x^*)\mu(x).
\]
Since $\mu$ is a contractive, positive functional, 
it follows from \cite[Proposition II.6.9.18]{blackadar:operator-algebras}
that $\mu(ax) = \mu(a)\mu(x)$ for every $a$ in $\C_0(\G)$. 
The same result, which is originally due to Choi \cite{choi:schwarz}, 
implies that  $\mu(xa) = \mu(x)\mu(a)$ because $\mu(x x^*)=\mu(x)\mu(x^*)$. 
\end{proof}

\begin{lemma} \label{lemma:mua}
Suppose that $X$ is symmetric. For every $\mu$ in $F_0$ and 
$a$ in $\C_0(\G)$ such that $\mu(a^* a)\ne 0$, the functional $\mu_a$ 
defined by
\[
\mu_a(b) = \frac{\mu(a^* b a)}{\mu(a^* a)}\qquad (b\in \C_0(\G))
\]
is in $F_0$.
\end{lemma}

\begin{proof}
The functional $\mu_a$ is clearly a state, so 
it is enough to show that $(\id\otimes \mu_a)\cop(x) = x$ 
for every $x$ in $X$. 
Suppose first that $a^* = (\id\otimes \omega)W$ for some 
$\omega$ in $\B(\ltwo(\G))_*$.
To simplify notation, put $\alpha = 1/\mu(a^* a)$. 
Let $\sigma\in \B(\ltwo(\G))_*$. Now 
\begin{align*}
\sigma(\id\otimes \mu_a)\cop(x)
 &= \alpha(\sigma\otimes \mu)\bigl((1\otimes a^*)\cop(x)(1\otimes a)\bigr)\\
 &= \alpha(\sigma\otimes \mu\otimes \omega)
        \bigl(W\sw{23}(\cop(x)\otimes 1)(1\otimes a\otimes 1)\bigr)\\
 &= \alpha(\mu\otimes \omega)
        \bigl(W((\sigma\otimes\id)\cop(x)\otimes 1)W^*W(a\otimes 1)\bigr).
\end{align*}
Write $x'=(\sigma\otimes\id)\cop(x)$ and note that $x'\in X$ 
by left invariance.

As noted in the proof of Theorem~\ref{thm:subalgebra->subgroup}, 
there exists $y$ in $X$ such that $\epsilon(y) = 1$.
By Lemma~\ref{lemma:basic-F-0}, $\mu.y = \mu$.
Moreover, $\omega = \tau.K$ for some $\tau$ in $\B(\ltwo(\G))_*$ and 
$K$ in $\B_0(\ltwo(\G))$. Inserting these into the preceding calculation gives
\[
\sigma(\id\otimes \mu_a)\cop(x) = 
\alpha(\mu\otimes \tau)\bigl((y\otimes K)W(x'\otimes 1)W^*W(a\otimes 1)\bigr).
\]
Since $X$ is symmetric, $W(x'\otimes 1)W^*\in\M(X\otimes \B_0(\ltwo(\G)))$,
so $(y\otimes K)W(x'\otimes 1)W^*$ is in $X\otimes \B_0(\ltwo(\G))$.
We replace this term by $z\otimes K'$ with $z$ in $X$ and $K'$
in $\B_0(\ltwo(\G))$ (although in reality $(y\otimes K)W(x'\otimes 1)W^*$ 
is only approximated in norm by sums of simple tensors $z\otimes K'$). 
We are now left with
\[
\alpha(\mu\otimes \tau)\bigl((z\otimes K')W(a\otimes 1)\bigr)
= \alpha\mu\bigl(z((\id\otimes \tau.K')W) a\bigr),
\]
which by Lemma~\ref{lemma:basic-F-0} is equal to
\begin{align*}
&\alpha\mu(z)\mu\bigl(((\id\otimes \tau.K')W) a\bigr)
=\alpha\epsilon(z)(\mu\otimes \tau)\bigl((1\otimes K')W(a\otimes 1)\bigr)\\
&\qquad
= \alpha(\mu\otimes \tau)
     \bigl((1\otimes(\epsilon\otimes \id)(z\otimes K'))W(a\otimes 1)\bigr)\\
&\qquad
= \alpha(\mu\otimes \tau)\bigl((1\otimes(\epsilon\otimes \id)
         ((y\otimes K)W(x'\otimes 1)W^*))W(a\otimes 1)\bigr),
\end{align*}
where at the final stage we replaced $z\otimes K'$ back to its 
true form. Noting that $\epsilon\otimes \id$ is multiplicative 
on $\M(\C_0(\G)\otimes \B_0(\ltwo(\G)))$  
and that $(\epsilon\otimes \id)W = (\epsilon\otimes \id)W^* = 1$
\cite[Theorem~3.1]{bedos-tuset:amenability},
we get 
\[
\sigma(\id\otimes \mu_a)\cop(x) = \alpha(\mu\otimes \tau)
    \bigl((1\otimes\epsilon(y)\epsilon(x')K)W(a\otimes 1)\bigr).
\]
But $\epsilon(y) = 1$ and $\epsilon(x') = \sigma(x)$, so we have
\[
\sigma(\id\otimes \mu_a)\cop(x) = 
\alpha \mu \bigl(((\id\otimes \tau.K)W )a\bigr)\sigma(x)
=\alpha \mu(a^* a)\sigma(x) = \sigma(x).
\]
This finishes the special case $a^* = (\id\otimes \omega)W$. 
The general case follows by approximating $a^*$ in norm.
\end{proof}

\begin{theorem} \label{thm:F=kerpiT}
Suppose that $\G$ is a co-amenable quantum group
and that $X$ is non-zero, symmetric, left invariant 
C*-subalgebra of\/ $\C_0(\G)$. 
A state $\mu$ of\/ $\C_0(\G)$ is in $F_0$ if and only if
its GNS-representation is $X$-trivial.
Moreover, if\/ $(\H,\pi)$ is the compact subgroup associated with $X$,
then $F_0 = \pi^*\bigl(\st(\C(\H))\bigr)$ 
where $\pi^*\col \C(\H)^*\to \C_0(\G)^*$ is 
the adjoint of the quotient map $\pi$.
\end{theorem}

\begin{proof}
Let $\mu\in F_0$. Applying the GNS-construction 
to $\mu$, we get a representation $\rho$ of $\C_0(\G)$ on 
a Hilbert space $H$ and a cyclic vector $\xi\in H$ such that 
$\mu = \pair{\rho(\cdot)\xi\mid\xi}$. 
We claim that $\rho$ is $X$-trivial. 
Let $x\in X$. Then it follows from Lemma~\ref{lemma:mua}
that for every $a\in\C_0(\G)$ 
\[
\pair{\rho(x)\rho(a)\xi\mid\rho(a)\xi} 
= \mu(a^*xa) = \mu(a^*a)\epsilon(x) 
= \pair{\rho(a)\xi\mid\rho(a)\xi} \epsilon(x).
\]
Since $\xi$ is cyclic, we have 
\[
\pair{\rho(x)\zeta\mid\zeta} = \pair{\epsilon(x)\zeta\mid\zeta} 
\]
for every $\zeta$ in $H$. It then follows from the polarisation identity
that $\rho(x) = \epsilon(x)1$, that is, $\rho$ is $X$-trivial. 
The converse is clear from Lemma~\ref{lemma:basic-F-0}.

By the first statement every $\mu$ in $F_0$ factors through 
$\pi\col \C_0(\G)\to \C(\H)$. 
Therefore $\mu$ is of the form $\nu\pi$ where $\nu$ is 
a state of $\C(\H)$. 
Conversely, if $\nu$ is a state of $\C(\H)$,
then $\nu\pi$ is a state of $\C_0(\G)$ and $\nu\pi = \epsilon$ 
on $X$ because $\pi$ is $X$-trivial. 
\end{proof}

\section{Compact quantum subgroups of a 
         co-amenable quantum group are co-amenable}
\label{sec:coamen}

The purpose of this section is to prove that 
a compact quantum subgroup of a co-amenable quantum group 
$\G$ is also co-amenable. In the special case when $\G$ is compact,
the result has been proved by Tomatsu 
(see \cite[Lemma~2.11]{tomatsu:coideal}).
When applied to the co-com\-mu\-ta\-tive case, the result says that
the quotient group of an amenable locally compact group by 
a normal, open subgroup is also amenable, which is of course 
well known from the classical theory. The section forms an independent 
part of the paper. 

Suppose that $\G$ is a co-amenable quantum group,
$\H$ is a compact quantum group 
and  $\pi\col \C_0(\G)\to\C(\H)$ is 
a surjective $*$-homomorphism such that 
$(\pi\otimes\pi)\cop = \cop_\H\pi$. 
In particular, $\H$ could be a compact quantum subgroup 
of $\G$ and $\pi$ the subgroup morphism composed with 
the natural map $\C^u(\H)\to\C(\H)$.

Recall that $W\in\M(\C_0(\G)\otimes  \C_0(\dual\G))$, and define 
\[
\tau\col \lone_*(\H)\to \C(\dual\G), \qquad 
\tau(\omega) = (\omega\pi\otimes \id)W.
\]
The following lemma is straightforward. 

\begin{lemma} \label{lemma:tau-isomorphism}
The map $\tau$ is an injective $*$-homomorphism from $\lone_*(\H)$ 
into $\C(\dual \G)$. 
\end{lemma}

By the definition of $\C^u_0(\dual\H)$,
there is a $*$-homomorphism 
$\rho\col \C^u_0(\dual\H)\to \C(\dual\G)$
such that $\tau = \rho \lambda_u$ where 
$\lambda_u\col \lone_*(\H)\to \C^u_0(\dual\H)$ is the natural embedding. 
A quantum group is said to be \emph{discrete} if its dual quantum group
is compact. A discrete quantum group is always co-amenable
(see for example \cite[Proposition~4.1]{runde:compact-discrete}),
so we may identify $\C^u_0(\dual\H)$ with $\C_0(\dual \H)$ 
\cite[Theorem~3.1]{bedos-tuset:amenability}. 
In this identification, $\lambda_u(\omega) = (\omega\otimes \id)U$ 
where $U\in \B(\ltwo(\H)\otimes \ltwo(\H))$ 
is the left multiplicative unitary of $\H$.
Therefore $\rho\col \C_0(\dual\H)\to \C(\dual\G)$ and 
$\rho\bigl((\omega\otimes \id)U\bigr) = (\omega\pi\otimes \id)W$
for every $\omega$ in $\lone_*(\H)$. It follows that 
\begin{equation} \label{eq:U-W}
(\id\otimes\rho)U = (\pi\otimes\id)W
\end{equation}
in $\M(\C(\H)\otimes \C_0(\dual\G))$. 
It should be pointed out that the map $\rho$ can always be defined on 
the universal C*-algebra level, using the universal co-representation
from \cite{kustermans:universal}.
We shall prove that it is injective when $\H$ is compact.

For the next result, denote the co-multiplications of $\dual \G$ and
$\dual \H$ by $\dual\cop$ and $\dual\cop_{\H}$, respectively.

\begin{theorem} \label{thm:rho}
The map $\rho\col \C_0(\dual \H)\to\C(\dual\G)$ is an injective 
non-degenerate $*$-homo\-morphism such that
$(\rho\otimes\rho)\dual{\cop}_{\H} = \dual{\cop}\rho$.
\end{theorem}

\begin{proof}
To show that $\rho$ is non-degenerate,
we adapt the argument used in 
\cite[Proposition~6.1]{maes-vandaele:compact}.
Write $\widetilde W=(\pi\otimes\id)W$.
Note that if $\omega\in\lone_*(\H)$ and $u\in(\dom S_\H)^*$, 
then the functional $u.\omega\col a\mapsto\omega(au)$ is also in $\lone_*(\H)$.
Therefore, for every $u$ in $(\dom S_\H)^*$ and $v$ in $\C_0(\dual\G)$,
\[
(\omega\otimes\id)\bigl(\widetilde W(u\otimes v)\bigr)
=\bigl((u.\omega\otimes\id)\widetilde W\bigr) v
\in \rho\bigl(\C_0(\dual \H)\bigr)\C_0(\dual \G).
\]
Since  $\widetilde W$ is a unitary in 
$\M\bigl(\C(\H)\otimes\C_0(\dual \G)\bigr)$, the set
$\widetilde W\bigl(\C(\H)\otimes\C_0(\dual \G)\bigr)$ is dense in 
$\C(\H)\otimes\C_0(\dual \G)$.
It follows that $\rho\bigl(\C_0(\dual \H)\bigr)\C_0(\dual \G)$
is dense in $\C_0(\dual\G)$,
so $\rho$ is non-degenerate.

We shall use the survey \cite{maes-vandaele:compact}
as the basis of our notation for compact quantum groups.  
Let $B_0$ be the subspace of $\C(\H)$ spanned by the matrix coefficients  
$u^i_{pq}$ of the irreducible unitary representations  
of $\H$ (here $i\in I$ and $1\le p,q \le n(i)$).
Then $B_0$ is a Hopf $*$-algebra, which is dense in $\C(\H)$. 
As shown in \cite[Lemma~8.1]{maes-vandaele:compact},
for every $i$ in $I$ and $0\le p,q\le n(i)$, 
there is $\omega_{pq}^i$ in $\lone(\H)$ such that
$\omega_{pq}^i(u_{pq}^i) = 1$ and 
$\omega_{pq}^i = 0$ at the other coefficients. 
The functional $\omega_{pq}^i$  is in fact in $\lone_*(\H)$
because $(\omega_{pq}^i)^* = \omega_{qp}^i$
\cite[p. 108]{maes-vandaele:compact}.
Let $\dual B _0$ be the subspace of $\lone_*(\H)$ spanned by the elements
$\omega_{pq}^i$. Then $\dual B_0$ is a $*$-subalgebra, which is
$*$-isomorphic with the algebraic direct sum of the 
finite-dimensional matrix algebras 
$M_{n(i)}$. Therefore $\dual B_0$ has a unique C*-completion:
the C*-algebraic direct sum of the matrix algebras $M_{n(i)}$.
The restriction of $\tau$ is a $*$-isomorphism from $\dual B_0$ onto
a $*$-subalgebra of $\C(\dual\G)$ by Lemma~\ref{lemma:tau-isomorphism}.
On the other hand, the restriction of $\lambda_u$ is a $*$-isomorphism from
$\dual B_0$ onto a dense $*$-subalgebra of $\C_0(\dual \H)$. 
It follows from the uniqueness of the C*-completion of $\dual B_0$ 
that $\rho\col \C_0(\dual \H) \to \C(\dual \G)$ is injective. 

It remains to show that $(\rho\otimes\rho)\dual{\cop}_{\H} = \dual{\cop}\rho$.
For every $\omega$ in $\lone_*(\H)$,
\begin{align*}
\dual{\cop}\rho\bigl((\omega\otimes\id)U\bigr)
&=\dual{\cop}\bigl((\omega\pi\otimes\id)W\bigr)
= \dual W^*\bigl(1\otimes (\omega\pi \otimes\id)W\bigr)\dual W\\
&= \Sigma W\bigl((\omega\pi \otimes\id)W\otimes 1\bigr)W^*\Sigma
= \Sigma\bigl((\omega\pi\otimes\id\otimes\id)
                     W\sw{23}W\sw{12}W\sw{23}^*\bigr)\Sigma\\
&= \Sigma\bigl((\omega\pi\otimes\id\otimes\id) W\sw{12}W\sw{13}\bigr)\Sigma
\end{align*}
by the pentagonal equation. Applying \eqref{eq:U-W}, we have
\begin{align*}
\dual{\cop}\rho\bigl((\omega\otimes\id)U\bigr)
& = (\omega\pi\otimes\id\otimes\id)   W\sw{13}W\sw{12}
= (\omega\otimes\id\otimes\id)(\id\otimes\rho\otimes\rho)U\sw{13}U\sw{12}\\
&=(\rho\otimes\rho)\dual{\cop}_{\H}\bigl((\omega\otimes\id)U\bigr).
\end{align*}
The claim follows because the elements $(\omega\otimes\id)U$ with
$\omega$ in $\lone_*(\H)$ are dense in~$\C_0(\dual \H)$.
\end{proof}

A quantum group $\G$ is said to be \emph{amenable}
if there exists a state $m$ of $\linfty(\G)$ such that 
\[
m(\omega\otimes\id)\cop(x) = \omega(1)m(x)\qquad
(\omega\in\lone(\G),\,x\in\linfty(\G)).
\]
Such a state is called a \emph{left invariant mean} on $\linfty(\G)$. 

As shown in \cite{bedos-tuset:amenability}, 
the dual quantum group of a co-amenable quantum group is amenable. 
The converse -- whether the dual quantum group of an amenable
quantum group is co-amenable -- is still open in general. 
The converse is known to be true in the commutative case
(i.e., for classical groups) and in the discrete case 
\cite{tomatsu:discrete-amenable}. We shall use the latter result in 
the proof of the following theorem.

\begin{theorem} \label{thm:co-amenable-subgroup}
Suppose that $\G$ is a co-amenable quantum group,
$\H$ is a compact quantum group and  $\pi\col \C_0(\G)\to\C(\H)$ is 
a surjective $*$-homomorphism such that $(\pi\otimes\pi)\cop = \cop_\H\pi$. 
Then also $\H$ is co-amenable. In particular, a compact quantum subgroup of 
a co-ame\-na\-ble  quantum group is co-amenable.
\end{theorem}

\begin{proof}
By \cite{tomatsu:discrete-amenable}, it suffices to show 
that $\dual \H$ is amenable. As already noted in the 
proof of Theorem~\ref{thm:rho}, $\C_0(\dual\H)$ is a C*-algebraic 
direct sum of full matrix algebras.
Therefore its multiplier algebra $\C(\dual\H)$
is the C*-algebraic direct product of the same algebras
and is $*$-isomorphic with the universal enveloping von Neumann algebra
$\C_0(\dual\H)^{**}$ of $\C_0(\dual\H)$.  
It follows from Theorem~\ref{thm:rho} that
the strict extension of $\rho$ from $\C(\dual \H)$ to $\C(\dual \G)$
is injective \cite{lance:hilbert-cstar}. 
On the other hand, $\rho$ extends to a normal $*$-homomorphism 
from $\C_0(\dual\H)^{**}$ onto the double commutant 
$\rho(\C_0(\dual\H))''\sub \B(\ltwo(\G))$.
Due to non-degeneracy, the normal extension coincides with the strict 
extension. It follows that the extension, which we still denote by $\rho$,
is  a $*$-isomorphism between the von Neumann algebras 
$\linfty(\dual \H) = \C(\dual \H)$ and $\rho(\dual \H)''$,
and so it is a homeomorphism with respect to the weak* topologies. 

Since $\G$ is co-amenable,
$\dual \G$ is amenable \cite{bedos-tuset:amenability}, 
and so there is a left invariant mean $m$ on $\linfty(\dual \G)$.
We shall show that $m\rho$ is a left invariant mean on 
$\linfty(\dual \H)$. Let $\omega\in\lone(\dual\H)$.
Since $\rho$ is a weak*-homeomorphism, there is $\sigma$ in $\lone(\dual \G)$
such that $\omega = \sigma\rho$. Now for every $x$ in $\linfty(\dual \H)$
\begin{align*}
m\rho(\omega\otimes \id)\dual{\cop}_{\H}(x) 
&= m (\sigma\otimes \id)(\rho\otimes \rho)\dual{\cop}_{\H}(x) \\
&=m (\sigma\otimes \id)\dual{\cop}\rho(x)
=\sigma(1) m\rho(x) 
=\omega(1) m\rho(x)
\end{align*}
because $m$ is left invariant. We conclude that $m\rho$ is a left 
invariant mean on $\linfty(\dual \H)$ and so $\dual\H$ is amenable. 
\end{proof}

\section{Left invariant C*-subalgebras from compact quantum subgroups}
\label{sec:subgroup->subalgebra}

In this section we proceed to the opposite direction from 
section~\ref{sec:subalgebra->subgroup}.
That is to say that we start with a compact quantum subgroup $(\H,\pi)$ 
of a co-amenable quantum group $\G$ 
and construct a left invariant C*-subalgebra of $\C_0(\G)$. 
In \cite{vaes:imprimitivity}, Vaes gave a general construction of a 
quantum homogeneous space associated with a closed quantum subgroup. 
Although related, that construction is much more elaborate than 
the construction given here: the reason is of course the more specialised 
setting of the present paper. 

Recall from the previous section that $\H$ is 
co-amenable and so $\pi\col \C_0(\G)\to\C(\H)$.
Let $F_0 = (\ker\pi)^\perp\cap \st(\C_0(\G))$. Then
\[
X = \set{x\in \C_0(\G)}{(\id\otimes \mu)\cop(x)=x\text{ for every }\mu\in F_0}.
\]
is a closed subspace of $\C_0(\G)$, and $X$ is also closed 
under involution. It is also easy to see that $X$ is left 
invariant. We shall show that $X$ is a 
C*-subalgebra of $\C_0(\G)$ by showing that it is closed under 
multiplication. 

Let $\phi_{\H}$ be the Haar state of $\H$. Put $\theta = \phi_{\H} \pi$,
and note that $\theta$ is in $F_0$. Define $P \col \C_0(\G) \to \C_0(\G)$ by 
\[
P(a) = (\id\otimes \theta)\cop (a)\qquad(a\in \C_0(\G)).
\]

\begin{lemma} \label{lemma:P}
\begin{enumerate}
\item The image of $P$ is $X$.
\item $P^2 = P$.
\item $P\bigl(P(a)P(b)\bigr) = P(a)P(b)$ for every $a$ and $b$ in $\C_0(\G)$.
\end{enumerate}
\end{lemma}

\begin{proof}
Each $\mu$ in $F_0$ is of the form $\mu = \mu' \pi$ 
where $\mu'$ is a state of $\C(\H)$. Since $\phi_{\H}$ is left invariant,  
\[
\mu*\theta = (\mu\otimes \theta)\cop = 
(\mu'\otimes \phi_{\H})(\pi\otimes \pi)\cop = 
(\mu'\otimes \phi_{\H})\cop_{\H}\pi = \phi_{\H}\pi = \theta.
\]
Then, for every $a$ in $\C_0(\G)$, 
\begin{align*}
(\id\otimes \mu)\cop(P(a)) 
&= (\id\otimes\mu\otimes\theta)(\cop\otimes \id)\cop (a)\\
&= (\id\otimes\mu *\theta)\cop (a)
= (\id\otimes\theta)\cop (a) = P(a).
\end{align*}
Therefore the image of $P$ is contained in 
$X$. Since $\theta$ is in $F_0$, it follows that $P(x) = x$
for every $x$ in $X$. This proves the first statement.

The second statement, that $P^2 = P$, follows from the identity 
$\theta*\theta = \theta$.

To prove the third statement, let $a ,b\in \C_0(\G)$. Then 
\begin{align*}
P\bigl(P(a)P(b)\bigr) 
&=(\id\otimes \theta)\cop\bigl((\id\otimes \theta)\cop (a)
                               (\id\otimes \theta)\cop (b)\bigr)\\
&=(\id\otimes\theta)\bigl(
            (\id\otimes\id\otimes\theta)(\id\otimes\cop)\cop (a)
(\id\otimes\id\otimes\theta)(\id\otimes\cop)\cop (b)\bigr)\\
&=(\id\otimes\phi_{\H})\bigl(
    (\id\otimes\id\otimes\phi_{\H})(\id\otimes(\pi\otimes\pi)\cop)\cop(a)\\
&\qquad\times
(\id\otimes\id\otimes\phi_{\H})(\id\otimes(\pi\otimes\pi)\cop)\cop(b)\bigr)\\
&=(\id\otimes\phi_{\H})\bigl(
    (\id\otimes\id\otimes\phi_{\H})(\id\otimes\cop_{\H}\pi)\cop(a) \\
&\qquad \times(\id\otimes\id\otimes\phi_{\H})(\id\otimes\cop_{\H}\pi)\cop(b)\bigr).
\end{align*}
Denote the unit map $\alpha\mapsto \alpha 1_{\H}\col \complex \to \C(\H)$
by $\eta$, and note that $(\id \otimes \phi_{\H})\cop_{\H} = \eta\phi_{\H}$.
Inserting this into the preceding calculation gives
\begin{align*}
P\bigl(P(a)P(b)\bigr) 
&=(\id\otimes\phi_{\H})\bigl((\id\otimes \eta\theta)\cop (a)
                          (\id\otimes \eta\theta)\cop (b)\bigr).\\
&=(\id\otimes\phi_{\H})\bigl(((\id\otimes \theta)\cop (a)
                           (\id\otimes \theta)\cop (b))\otimes 1_{\H}\bigr).\\
&=P(a)P(b)\phi_{\H}(1_{\H}) = P(a)P(b).
\end{align*}
\end{proof}

Recall that a map from a C*-algebra onto its C*-subalgebra 
is called a \emph{conditional expectation} if 
it is a projection of norm $1$. 
This is not the traditional definition but equivalent to it 
(see \cite[subsection~II.6.10]{blackadar:operator-algebras}). 
Every conditional expectation $E$ is completely positive 
and satisfies $E(aE(b)) = E(a)E(b) = E(E(a)b)$ for every $a$ and~$b$. 

\begin{theorem} \label{thm:subgroup->subalgebra}
Suppose that $(\H,\pi)$ is a compact quantum subgroup of 
a co-ame\-na\-ble  quantum group $\G$. 
Let $X$ be the subspace associated with $\H$.
Then $X$ is a non-zero, symmetric, left invariant C*-subalgebra of\/ $\C_0(\G)$.
Moreover, the map $P$ is a conditional expectation from $\C_0(\G)$ onto $X$ such that
$(\id \otimes P)\cop = \cop P$.
\end{theorem}

\begin{proof}
We have already noticed that $X$ is a closed, left invariant 
subspace of $\C_0(\G)$ which is closed under involution. 
It follows from the third statement 
of the preceding lemma that $X$ is closed under multiplication and 
is therefore a C*-subalgebra of~$\C_0(\G)$. The C*-subalgebra $X$ is non-zero
because if $a\in \C_0(\G)$ such that $\phi_{\H}\pi(a)\ne 0$, then $P(a) \ne 0$.

The map $P\col \C_0(\G)\to X$ is a conditional expectation
because it is a surjective projection of norm $1$
by Lemma~\ref{lemma:P}.
If $(e_i)$ is a bounded approximate identity in $\C_0(\G)$,
then $(P(e_i))$ is a bounded approximate identity in $X$. 
It follows that the completely positive map
$\id\otimes P\col \C_0(\G)\otimes \C_0(\G)\to \C_0(\G)\otimes X$
is non-degenerate, so it has a strict extention to a 
map $\id \otimes P\col \M(\C_0(\G)\otimes \C_0(\G))\to \M(\C_0(\G)\otimes X)$.
The identity $(\id \otimes P)\cop = \cop P$ follows immediately 
from the definition of $P$ and the co-associativity of $\cop$.

It remains to show that $X$ is symmetric. 
Note first that the strict extension 
$P\otimes \id\col \M\bigl(\C_0(\G)\otimes \B_0(\ltwo(\G))\bigr)\to 
\M\bigl(X\otimes \B_0(\ltwo(\G))\bigr)$
is also a conditional expectation.
Therefore it suffices to show that
\begin{equation*} 
(P\otimes \id)(W(x\otimes 1)W^*) = W(x\otimes 1)W^*
\end{equation*}
for every $x$ in $X$. By the pentagonal equation, 
\begin{align*}
&(P\otimes \id)\bigl(W(x\otimes 1)W^*\bigr) 
=(\id\otimes\theta\otimes \id)\bigl(W_{12}^* W_{23}(1\otimes x\otimes 1)
                                   W_{23}^* W_{12}\bigr)\\
&\qquad 
 =(\id\otimes\theta\otimes \id)\bigl(W_{13} W_{23}W_{12}^*(1\otimes x\otimes 1)
                                   W_{12}W_{23}^* W_{13}^*\bigr)\\
&\qquad= W\bigl((\id\otimes\theta\otimes \id)
         (W_{23}(\cop(x)\otimes 1)W_{23}^*) \bigr)W^*,
\end{align*}
so all we need to show is that
\[
(\id\otimes\theta\otimes \id)(W_{23}(\cop(x)\otimes 1) W_{23}^* )
=x\otimes 1.
\]
Take an arbitrary element in $\B(\ltwo(\G))_*$ and factorise it
as $K_2.\omega.K_1$ with $\omega\in \B(\ltwo(\G))_*$ and 
$K_1, K_2\in \B_0(\ltwo(\G))$. 
Pick $a$ in $\C_0(\G)$ such that $\pi(a) = 1_{\H}$ and 
note that, as $\theta = \phi_{\H}\pi$, we have $\theta = a.\theta.a$. Now
\begin{align*}
&(\id\otimes\theta\otimes K_2.\omega.K_1)
  (W_{23}(\cop(x)\otimes 1)W_{23}^* )\\
&\qquad
= (\id\otimes\theta\otimes \omega)
  \bigl((1\otimes a\otimes K_1)W_{23}(\cop(x)\otimes 1)W_{23}^*
        (1\otimes a\otimes K_2)\bigr).
\end{align*}
Since $W\in \M\bigl(\C_0(\G)\otimes \B_0(\ltwo(\G))\bigr)$, we may replace 
$(1\otimes a\otimes K_1)W_ {23}$ by $1\otimes b\otimes K_3$
and $W_{23}^*(1\otimes a\otimes K_2)$ by $1\otimes c\otimes K_4$
where $b,c\in \C_0(\G)$ and $K_3, K_4\in\B_0(\ltwo(\G))$. 
Then we are left with the term 
\begin{align*}
&(\id\otimes\theta\otimes \omega)
  \bigl((1\otimes b\otimes K_3)(\cop(x)\otimes 1)
                                  (1\otimes c\otimes K_4)\bigr)\\
&\qquad= \omega(K_3K_4)(\id\otimes\phi_{\H}\pi)
        \bigl((1\otimes b)\cop(x)(1\otimes c)\bigr)\\
&\qquad= \omega(K_3K_4)(\id\otimes\phi_{\H})\bigl((1\otimes \pi(b))
         (\id\otimes\pi)\cop(x)(1\otimes \pi(c))\bigr).
\end{align*}
Now $\pi P = \eta\theta = \eta\epsilon P$ where $\eta\col \complex\to \C(\H)$ 
is the unit map and $\epsilon\col \C_0(\G)\to \complex$ the co-unit. Therefore 
\[
(\id\otimes\pi)\cop(x) 
= (\id\otimes \pi P)\cop(x)
= (\id\otimes \epsilon P)\cop(x)\otimes 1_{\H}
= x\otimes 1_{\H}.
\]
Applying this to the ongoing calculation gives
\begin{align*}
&(\id\otimes\theta\otimes K_2.\omega.K_1)(W_{23}(\cop(x)\otimes 1)W_{23}^* )
=\omega(K_3K_4)(\id\otimes\phi_{\H})(x\otimes \pi(bc))\\
&\quad= (\phi_{\H}\pi\otimes \omega)(bc\otimes K_3 K_4)x
       = (\theta\otimes \omega)\bigl((a\otimes K_1)W W^*(a\otimes K_2)\bigr)x
       = \omega(K_1K_2)x.
\end{align*}
Since the functional $K_2.\omega.K_1$ in $\B(\ltwo(\G))_*$ is arbitrary, 
it follows that
\[
(\id\otimes\theta\otimes \id)(W_{23}(\cop(x)\otimes 1) W_{23}^* )
 = x\otimes 1.
\]
as required.
\end{proof}

\section{Uniqueness results}

We have seen that every non-zero, left invariant C*-subalgebra 
of a co-amenable quantum group gives rise to a compact quantum subgroup 
and vice versa. In this section we show uniqueness results related to
these two constructions. 

Let $\G$ be a co-amenable quantum group.
For every non-zero, left invariant C*-subalgebra $X$ of $\C_0(\G)$,
denote the compact quantum subgroup of $\G$ associated with $X$ 
by $(\H_X,\pi_X)$ (see Theorem~\ref{thm:subalgebra->subgroup}). 
Conversely, for every compact quantum subgroup $\H$ 
of $\G$, denote the left invariant C*-subalgebra of $\C_0(\G)$ associated with 
$\H$ by $X_{\H}$ (see Theorem~\ref{thm:subgroup->subalgebra}). 

It follows from Theorem~\ref{thm:subgroup->subalgebra}
that if $X_{\H_X} = X$, then $X$ is symmetric and 
there is a conditional expectation 
$P\col \C_0(\G)\to X$ such that $(\id\otimes P)\cop = \cop P$. 
So the existence of such a conditional expectation is 
a consequence of any uniqueness result. The existence is assumed in 
the following result, leaving some room for improvement. 
(As we shall see in section~\ref{sec:co-commutative}, 
the symmetry condition is absolutely necessary.)

\begin{theorem} \label{thm:1st-unique}
Suppose that $\G$ is a co-amenable quantum group and that
$X$ is a non-zero, symmetric, left invariant C*-sub\-al\-gebra of\/ $\C_0(\G)$.
If there is a conditional expectation $P\col \C_0(\G)\to X$ such that 
$(\id\otimes P)\cop = \cop P$, then $X_{\H_X} = X$.
\end{theorem}

\begin{proof}
By definition,
\[
X_{\H_X} = \set{x\in \C_0(\G)}{(\id\otimes \mu)\cop(x) = x
                \text{ for every }\mu\in(\ker \pi_{X})^\perp\cap \st(\C_0(\G))}.
\]
Following section~\ref{sec:subalgebra->subgroup}, let 
\[
F_0 = \set{\mu\in \st(\C_0(\G))}
        {(\id\otimes \mu)\cop(x)= x\text{ for every }x\in X}.
\]
By Theorem~\ref{thm:F=kerpiT}, $F_0 = (\ker \pi_{X})^\perp\cap \st(\C_0(\G))$
and it follows that 
\[
X_{\H_X} = \set{a\in \C_0(\G)}{(\id\otimes \mu)\cop(a) = a
                               \text{ for every }\mu\in F_0}.
\]
Therefore $X\sub X_{\H_X}$. Conversely, let $a\in X_{\H_X}$. 
Since $\epsilon P \in F_0$, we have
\[
a = (\id\otimes \epsilon P)\cop(a) = (\id\otimes\epsilon)\cop P(a) = P(a),
\]
and so $a\in X$.
\end{proof}

\begin{theorem} \label{thm:2nd-unique}
Suppose that $(\H,\pi)$ is a compact quantum subgroup of 
a co-ame\-na\-ble quantum group $\G$. 
Then $(\H_{X_{\H}},\pi_{X_{\H}}) \cong (\H,\pi)$.
\end{theorem}

\begin{proof}
Denote the Haar state of $\H$ by $\phi_{\H}$, and write
$\theta_{\H} = \phi_{\H}\pi$ and $P_{\H} = (\id\otimes \theta_{\H})\cop$.
By Theorem~\ref{thm:subgroup->subalgebra},
$P_{\H}$ is a conditional expectation from $\C_0(\G)$ onto $X_{\H}$ 
such that $(\id\otimes P_{\H})\cop = \cop P_{\H}$.
Moreover, $\psi P_{\H} = \psi$ where $\psi$ is the right Haar 
weight of~$\C_0(\G)$. Let $\phi_{X_\H}$ be the Haar state of $\H_{X_{\H}}$ and 
put $\theta_{X_\H} = \phi_{X_\H}\pi_{X_{\H}}$. Now
$P_{X_\H} = (\id\otimes \theta_{X_\H})\cop$ is a 
conditional expectation from $\C_0(\G)$ onto $X_{\H_{X_{\H}}}$ 
such that $\psi P_{X_\H} = \psi$.
By Theorem~\ref{thm:1st-unique}, $X_{\H_{X_{\H}}} = X_{\H}$,
so both $P_{\H}$ and $P_{X_\H}$ are conditional expectations 
from $\C_0(\G)$ onto $X_{\H}$. Since a conditional expectation $P$ from 
$\C_0(\G)$ onto $X_{\H}$ such  that $\psi P = \psi$ is unique 
\cite[Proposition~II.6.10.10]{blackadar:operator-algebras},
it follows that $P_{X_\H} = P_{\H}$. Therefore 
$\theta_{X_\H} = \epsilon P_{X_\H} = \epsilon P_{\H} = \theta_{\H}$. 

Now let $a\in \ker \pi$. Then 
\[
0 = \phi_{\H}\bigl(\pi(a)^* \pi(a) \bigr) 
  =\theta_{\H}(a^* a ) = \theta_{X_\H}(a^*a) 
  = \phi_{X_\H}\bigl(\pi_{X_{\H}}(a)^*\pi_{X_{\H}}(a)\bigr). 
\]
Since $\phi_{X_\H}$ is faithful, we must have $\pi_{X_\H}(a) = 0$.
Therefore,  $\ker\pi\sub\ker\pi_{X_{\H}}$. 

Conversely, let $F_0 = \set{\mu\in \st(\C_0(\G))}{(\id\otimes \mu)\cop(x)= 
x\text{ for every }x\in X_{\H}}$, and recall that 
$\ker\pi_{X_{\H}} = \set{a\in\C_0(\G)}{\mu(a)=0\text{ for every }\mu\in F_0}$ 
by Theorem~\ref{thm:F=kerpiT}.
If $\nu$ is a state of $\C(\H)$, 
then $\nu \pi P_{\H} = \nu(1)\theta_{\H} = \theta_{\H}$.
It follows that $(\id\otimes \nu\pi)\cop P_{\H} = P_{\H}$, and so 
$\nu\pi\in F_0$.
We conclude that $\ker\pi_{X_{\H}} \sub \ker\pi$. 

Since $\ker \pi_{X_{\H}} = \ker \pi$, there is a $*$-isomorphism 
$\rho\col \C(\H)\to \C(\H_{X_{\H}})$ such that $\rho\pi = \pi_{X_{\H}}$.
In other words, $(\H_{X_{\H}},\pi_{X_{\H}}) \cong (\H,\pi)$.
\end{proof}

Putting the previous theorems together we obtain the following
correspondence result. It is similar to 
Theorem~3.18 of \cite{tomatsu:coideal}, which concerns 
co-amenable \emph{compact} quantum groups.

\begin{theorem}
Suppose that $\G$ is a co-amenable quantum group.
There is a one-to-one correspondence between
compact quantum subgroups of\/ $\G$ and 
non-zero, symmetric, left invariant C*-sub\-al\-gebras $X$ 
of\/ $\C_0(\G)$ with a conditional expectation $P\col \C_0(\G)\to X$ such that
$(\id\otimes P)\cop = \cop P$. 
\end{theorem}

\section{The co-commutative case} \label{sec:co-commutative}

Consider the case of a co-amenable, co-commutative quantum group,
that is, the dual of an \emph{amenable} locally compact group $G$.
As $G$ is amenable, the universal group C*-algebra $\cstar(G)$ 
is isomorphic with the reduced group C*-algebra $\cstarred(G)$. 
We shall use the simpler notation $\cstar(G)$. The dual space of 
$\cstar(G)$ is the Fourier--Stieltjes algebra $\B(G)$ of $G$, which 
consists of all matrix coefficients of unitary representations of $G$. 
The multiplication of $\B(G)$ induced by the co-multiplication 
on $\cstar(G)$ is the pointwise multiplication of functions. 
The Fourier algebra $\A(G)$, which consists of the coefficients 
of the left regular representation $\lambda$, is a closed ideal in $\B(G)$.
The dual space of $\A(G)$ is the von Neumann algebra $\vn(G)$
generated by $\cstar(G)$ in $\B(\ltwo(G))$. Eymard, who introduced the 
Fourier and the Fourier--Stieltjes algebras in \cite{eymard:fourier},  
also studied the action of $\B(G)$ on $\vn(G)$ defined by
\[
\pair{ua, v} = \pair{a,uv} \qquad(a\in\vn(G),\,u\in \B(G),\,v\in \A(G)).
\]
Restricted to $\cstar(G)$, this action coincides 
with the action $(u,a)\mapsto L_u(a)$ defined in section~\ref{sec:prel},
but we shall use the more conventional notation $ua$ instead of $L_u(a)$.
When $a = \lambda(f)$ for some $f$ in $\lone(G)$,
the action is given by pointwise multiplication 
of $\lone$-functions by continuous functions:
\[
u\lambda(f) = \lambda(uf)\qquad(u\in \B(G),\, f\in\lone(G)).
\]

Eymard also introduced the notion of \emph{support} of an operator $a$ 
in $\vn(G)$, denoted by $\supp a$ \cite{eymard:fourier}. 
A point $s$ in $G$ is in $\supp a$ if and only if, 
for every neighbourhood $U$ of $s$, there exists $u$ in $\A(G)$ 
such that $\supp u\sub U$ and $\pair{a,u}\ne 0$.

Suppose that $H$ is an \emph{open} subgroup of $G$. 
Then we may consider $\lone(H)$ as a subspace of $\lone(G)$,
and the group C*-algebra $\cstar(H)$ of $H$ may be identified with 
\[
X = \set{x\in \cstar(G)}{\supp x\sub H}  
  = \normcl{\lambda\bigl(\lone(H)\bigr)}.
\]
Note that $X$ is an invariant C*-subalgebra of $\cstar(G)$.
(Since the co-multiplication of $\cstar(G)$ is co-commutative,
$(\id\otimes u)\cop(a) = (u\otimes \id)\cop(a)$ 
for every $u$ in $\B(G)$ and $a$ in $\cstar(G)$, and so
one-sided invariance implies two-sided invariance.)
The following theorem shows that all non-zero, invariant
C*-subalgebras of $\cstar(G)$ are of this form. 
It is the dual version of the Lau--Losert theorem stated in the 
introduction. 

\begin{theorem} \label{thm:dual-of-LL}
Suppose that $G$ is an amenable locally compact group. There is a one-to-one 
correspondence between open subgroups $H$ of $G$ and non-zero, invariant 
C*-subalgebras $X$ of\/ $\cstar(G)$. The correspondence is given by
\begin{align}
X &= \set{x\in \cstar(G)}{\supp x\sub H} \label{eq:X}\\
H &= \bigcup_{x\in X} \supp x. \label{eq:H}
\end{align}
\end{theorem}

\begin{proof}
As already noted, an open subgroup $H$ determines a non-zero, 
invariant C*-subalgebra $X$ of $\cstar(G)$ via \eqref{eq:X}. 
Obviously $\bigcup_{x\in X} \supp x \sub H$ and the converse 
follows easily because $H$ is open. So we recover $H$ via \eqref{eq:H}.

Suppose now that we are given a non-zero, invariant C*-subalgebra $X$ 
of $\cstar(G)$. The double commutant $X''$ of $X$ is the weak* closure of 
$X$ in $\B(\ltwo(G))$. Put
\[
H = \set{s\in G}{\lambda(s)\in X''}
\]
and  
\[
Y = \set{x\in \cstar(G)}{\supp x\sub H}.
\]
It is easy to see that $H$ is a closed subgroup of $G$.
We show next that $X = Y$.

If $x\in X$ and $s\in \supp x$, then there is a 
net $(u_i)$ in $\A(G)$ such that $u_i x \to \lambda(s)$
in the weak* topology (by \cite[Proposition~4.4]{eymard:fourier}). 
Since $X$ is invariant, each  $u_i x\in X$. 
It follows that $s\in H$, and so $\supp x\sub H$. Therefore $X\sub Y$.

Conversely, let $y\in Y$. 
Since $\supp y \sub H$, it follows from 
\cite[Theorem~3]{takesaki-tatsuuma:duality-subgroups2}
that $y \in \lambda(H)''$. But $\lambda(H)''\sub X''$, so $y\in X''$. 
Suppose for awhile that the support of $y$ is compact.
Pick compactly supported $u$ in $\A(G)$ such that
$u = 1$ on a neighbourhood of $\supp y$. If $(y_i)$ is a net 
in $X$ converging to $y$ in the weak* topology, then
$uy_i \to uy = y$. By invariance,  the net $(uy_i)$ is in $X$, and
since $u$ is supported by some compact set $K$, each $uy_i$ is also 
supported by $K$. Let $v\in \A(G)$ such that $v=1$ on $K$.
Then for every $w$ in $\B(G)$
\[ 
\pair{uy_i, w} = \pair{uy_i, vw} \to \pair{y,vw} = \pair{y,w}.
\] 
In other words, $uy_i\to y$ in the weak topology of $\cstar(G)$. 
By the standard convexity argument, $y$ is in the norm closure
of the convex hull of $(uy_i)$, so $y\in X$. 
This solves the case of compactly supported $y$. 

Suppose now that $y$ in $Y$ is arbitrary. 
Since $G$ is amenable, there is a bounded approximate identity 
$(u_i)$ in $\A(G)$ such that each $u_i$ is compactly supported. 
Then the support of each $u_i y$ is compact and contained in $H$. 
It is well known that every element in $\cstar(G)$ 
has a factorisation $vz$ where $v\in \A(G)$ and $z\in\cstar(G)$: 
for compactly supported ones this is obvious and 
approximation gives the general case because $\A(G)\cstar(G)$ is closed
as a consequence of Cohen's factorisation theorem.
It follows that $u_i y\to y$ in norm. By the previous paragraph, 
each $u_i y$ is in $X$ and hence $y\in X$.

To show that $H$ is open, we use an argument inspired by 
\cite[proof of Lemma~3.2]{kaniuth-lau:separation-property}.
Assume that $H$ is not open. Let $x\in X$ and $\epsilon>0$. 
Choose a compactly supported continuous function $f$ on $G$ such that 
$\|x-\lambda(f)\| < \epsilon$. Since we assume that $H$ is not open, 
the left Haar measure $|H|$ of $H$ is $0$.
(If $|H|\ne 0$, there is a compact set $K\sub H$ with $|K|>0$, and 
so $g := 1_K*1_{K\inv}$ is a non-zero, continuous function supported by $H$.
Therefore $g\inv((0,\infty))$ is a non-empty open set contained
in $H$ and hence $H$ is open.)
So there is an open set $U$ such that $|U|<\epsilon/\|f\|_{\text{sup}}$
and $\supp f\cap H\sub U$. Now put $g = f 1_{G\sm U}$, 
where $1_{G\sm U}$ denotes the characteristic function of $G\sm U$. 
Then 
\[
\|x  - \lambda(g)\| < \epsilon + \|\lambda(f)-\lambda(g)\|
\le \epsilon + \|f-g\|_1 \le \epsilon + \|f\|_{\text{sup}} |U| < 2\epsilon.
\]

Let
\[
I(H) = \set{u\in \A(G)}{u = 0\text{ on }H}.
\]
Now $I(H)^\perp=\lambda(H)''$
by \cite[section~3]{kaniuth-lau:separation-property},
and $\lambda(H)''= X''$ 
by \cite[Theorem~3]{takesaki-tatsuuma:duality-subgroups2}.
It then follows from \cite[Theorem~1.3]{for-kan-lau-spr} that
there is a (completely) bounded projection $P\col \vn(G)\to X''$ 
such that 
\[
P(u a) = u P(a)\qquad(u\in \A(G), a\in \vn(G)).
\]
Since $\supp f\cap (G\sm U)$ is compact and $H$ is closed,
there is $u$ in $\A(G)$ such that $u = 1$ on $\supp f\cap (G\sm U)$ 
and $u = 0$ on $H$. For every $v$ in $\A(G)$, the function $uv$ 
vanishes on $H$, so
\[
0 =  \pair{P(\lambda(g)), uv}= \pair{u P(\lambda(g)), v}
  = \pair{P(\lambda(ug)),v } = \pair{P(\lambda(g)), v}.
\]
Hence $P(\lambda(g)) = 0$, and so 
\[
\|x\| = \|P(x) - P(\lambda(g))\| \le 2\|P\|\epsilon.
\] 
Since $\epsilon>0$ is arbitrary, $x = 0$.
Therefore $X = \{0\}$, which is in contradiction with the hypothesis
of the theorem. We conclude that $H$ is open. 

So the characterisation \eqref{eq:X} of $X$ holds for some open subgroup
$H$, and the openness of $H$ guarantees that, in fact, $H$ 
is of the form given in \eqref{eq:H}. 
\end{proof}

In the case of a locally compact \emph{abelian} group $G$, 
the subgroup duality $H\mapsto H^\perp$ maps compact subgroups of $G$ to 
open subgroups of $\dual G$ and vice versa. 
The change of compact subgroups from the Lau--Losert theorem 
to open subgroups in the preceding theorem reflects this duality 
between open and compact subgroups. 

We now fix a non-zero, invariant C*-subalgebra $X$ of $\cstar(G)$.
By the preceding theorem $X = \cstar(H)$ for some open subgroup $H$ of $G$.
Let $K$ be the conjugate closure of $H$, i.e., the smallest normal
subgroup of $G$ containing $H$. Note that also $K$ is open.
It can be shown that a non-degenerate $*$-homomorphism 
$\rho\col \cstar(G)\to \M(A)$ is $X$-trivial 
if and only if $\rho = 1$ on $K$.
It follows that in this case
the construction of section~\ref{sec:subalgebra->subgroup}
leads to the compact quantum subgroup $\H = (G/K)\latedual$
(so that $\C(\H)\cong \cstar(G/K)$).
This example shows that different invariant 
C*-subalgebras can lead to the same
compact quantum subgroup.

We move on to characterise the case when $X$ is symmetric. 
As is to be expected, this is the case exactly when
$H$ is normal.
Following the notation of section~\ref{sec:subalgebra->subgroup}, let 
\[
F_0 = \set{u\in\posdef_0(G)}{ux = x\text{ for every $x$ in $X$}}
\]
where $\posdef_0(G)$ denotes the state space of $\cstar(G)$.
Note that $\posdef_0(G)$  consists of 
all continuous positive definite functions 
on $G$ with value $1$ at the identity $e$ of $G$.

\begin{lemma} \label{lemma:1onH}
A function $u$ in $\posdef_0(G)$ is in $F_0$ 
if and only if $u = 1$ on~$H$.
\end{lemma}

\begin{proof}
Suppose first that $u\in F_0$. Let $h\in H$, and $f$ be 
a continuous function such that $f(h) = 1$
and the support of $f$ is a compact subset of $H$. 
Since $\lambda(f)\in X$, we have 
$\lambda(f) = u\lambda(f) = \lambda(uf)$.
Since both $f$ and $uf$ are continuous, 
$uf = f$ and so $u(h) = 1$.

Conversely, suppose that $u\in\posdef_0(G)$ such that $u = 1$ on $H$. 
Let $x \in X$. Then $u-1 = 0$ on a neighbourhood of $\supp x$
(namely, on $H$), so $(u-1)x = 0$ by \cite[Proposition 4.8]{eymard:fourier}. 
\end{proof}

\begin{lemma}
Every continuous positive definite function $u$ on $G$ with $u = 1$ on $H$
is constant on all left and right cosets of $H$. 
\end{lemma}

\begin{proof}
There is a unitary representation $\rho$ of $G$ on a Hilbert space $K$ 
and a unit vector $\xi$ in $K$ such that 
$u(s) = \pair{\rho(s)\xi\mid\xi}$ for every $s$ in $G$. 
For every $h$ in $H$, 
\[ 
1 = u(h) = \pair{\rho(h)\xi\mid\xi} \le \|\rho(h)\xi\| \le  1,
 \]
which implies that $\rho(h)\xi = \xi$. Therefore 
$u(sh) = \pair{\rho(sh)\xi\mid\xi} = \pair{\rho(s)\xi\mid\xi} = u(s)$
and $u(hs) = u(s)$ for every $s$ in $G$ and $h$ in~$H$.
\end{proof}

\begin{lemma} \label{lemma:const-on-cosets}
A function $u$ in $\B(G)$ is in $F = \linsp F_0$ 
if and only if $u$ is constant on all left and right cosets of $H$.
Moreover, $F$ is weak*-closed in $\B(G)$.
\end{lemma}

\begin{proof}
It follows from the preceding two lemmas that
every $u$ in $F$ is constant on all left and right cosets of $H$. 

Conversely, let $u\in \B(G)$ such that $u$ 
is constant on all left and right cosets of $H$.
Put $\tilde u(s) = \conj{u(s\inv)}$ for every $s$ in $G$.
Then $u = u_1 + i u_2$ where $u_1 = (u + \tilde u)/2$
and $u_2 = (u-\tilde u)/2i$. 
Since $u$ is constant on both the left and the right cosets 
of $H$, so are $u_1$ and $u_2$. It follows that we may assume without 
loss of generality that $u = \tilde u$; that is, $u$ is hermitian.
Let $u = u^+ - u^-$ be the Jordan decomposition of $u$,
so both $u^+$ and $u^-$ are positive definite. 
By \cite[Lemme 2.12]{eymard:fourier}, for any given $\epsilon>0$,
there are $\alpha_k$ in $\complex$ and $s_k$ in $G$ such that 
\[
\biggl\| u^+(t) - \sum_{k=1}^n \alpha_k u(t s_k)\biggr\| < \epsilon
\qquad\text{for every $t$ in $G$}.
\]
Now for every $h$ in $H$,
\[
| u^+(h) - u^+(e)| \le 2\epsilon
+ \biggl|\sum_{k=1}^n \alpha_k u(hs_k) - 
         \sum_{k=1}^n \alpha_k u(s_k)\biggr| = 2\epsilon.
\]
It follows that $u^+(h) = u^+(e)$, so either $u^+ = 0$
or $u^+/u^+(e)$ is in $F_0$.
The same holds for $u^-$ and hence $u\in F$.

To show that $F$ is weak*-closed, it is enough to 
consider bounded nets $u_i\to u$ where 
each $u_i$ is a hermitian element in $F$. 
For each $i$, let $u_i = u_i^+-u_i^-$ be the Jordan
decomposition. Since the nets $(u_i^+)$ and $(u_i^-)$ are bounded,
a subnet argument shows that there are positive 
elements $v^+$ and $v^-$ in $\B(G)$ such that $u_i^+\to v^+$, 
$u_i^-\to v^-$ and $v^+ - v^- = u$. 
By the beginning of the proof, the nets $(u_i^+)$ and 
$(u_i^-)$ are in $\R^+F_0$. Since $\R^+F_0$ is 
weak*-closed, $u\in F$ as required. 
\end{proof}

For the following theorem, define 
\[
F_\perp = \set{a\in\cstar(G)}{\pair{u,a}=0 \text{ for every $u$ in $F
    = \linsp F_0$}}.
\]

\begin{theorem} \label{thm:dual-normal}
Suppose that $H$ is an open subgroup of an amenable
locally compact group $G$. Then the following statements are equivalent:
\begin{enumerate}
\item \label{item:H-normal} $H$ is normal,
\item \label{item:symmetric} $\cstar(H)$ is symmetric,
\item \label{item:Fperp-ideal} $F_\perp$ is an ideal in $\cstar(G)$.
\end{enumerate}
When $H$ is normal, the compact quantum subgroup $(\H,\pi)$ 
of\/ $\G=\dual G$ arising from $X = \cstar(H)$ 
is isomorphic with the compact quantum subgroup 
$(G/H)\latedual$ of $\dual G$.
\end{theorem}

\begin{proof}
To show that (\ref{item:H-normal}) implies (\ref{item:symmetric}),
suppose that $H$ is normal. Recall that the left multiplicative 
unitary associated with $\cstar(G)$ is given by 
\[
\dual W\xi(s,t) = \xi(ts, t)\qquad(\xi\in\ltwo(G\times G), s,t\in G).
\]
We should show that $\dual W(x\otimes 1)\dual{W}^*$ is
in $\M\bigl(\cstar(H)\otimes \B_0(\ltwo(G))\bigr)$ 
for every $x$ in $\cstar(H)$. 
Now $\dual W$ is in $\M\bigl(\cstar(G)\otimes \C_0(G)\bigr)$
and so is $\dual W(x\otimes 1)\dual{W}^*$.
We may identify $\M\bigl(\cstar(G)\otimes \C_0(G)\bigr)$ with
$\C^{\text{str}}\bigl(G,\M(\cstar(G))\bigr)$, 
the bounded functions from $G$ to $\M(\cstar(G))$ that 
are continuous with respect to the strict topology on $\M(\cstar(G))$:
an element $T$ in $\M\bigl(\cstar(G)\otimes \C_0(G)\bigr)$ 
is identified with the function
\[
s\mapsto (\id\otimes\delta_s)T\col G\to \M(\cstar(G)),
\]
where $\delta_s$ is the point mass at $s\in G$. 
If we can show that $\dual W(x\otimes 1)\dual{W}^*$ is in 
fact in $\C^{\text{str}}\bigl(G,\M(\cstar(H))\bigr)$, 
then the claim follows. 
Since $\id\otimes \delta_s$ is multiplicative on 
$\M\bigl(\cstar(G)\otimes \C_0(G)\bigr)$ and 
$(\id\otimes \delta_s)\dual W = \lambda(s\inv)$,
we have
\[
(\id\otimes \delta_s)\bigl(\dual W(x\otimes 1)\dual{W}^*\bigr)
= \lambda(s\inv) x \lambda(s).
\]
But $\supp(\lambda(s\inv) x \lambda(s))\sub s\inv \supp(x) s$
(by \cite[Proposition 4.8]{eymard:fourier}), 
so the normality of $H$ implies that $\lambda(s\inv) x \lambda(s)\in\cstar(H)$.
It follows that $\dual W(x\otimes 1)\dual{W}^*$ is in 
$\C^{\text{str}}\bigl(G,\M(\cstar(H))\bigr)$, and so 
$\cstar(H)$ is symmetric.

If $\cstar(H)$ is symmetric, then $F_\perp = \ker\pi$ 
by Theorem~\ref{thm:F=kerpiT}.
So (\ref{item:symmetric}) implies (\ref{item:Fperp-ideal}).

To show that  (\ref{item:Fperp-ideal}) implies (\ref{item:H-normal}), 
suppose that $F_\perp$ is an ideal. 
Since $H$ is an open subgroup, the characteristic function
$1_H$ is a continuous positive definite function 
\cite[subsection 32.43]{hewitt-ross:vol2}. By Lemma~\ref{lemma:1onH}, 
$1_H\in F$. Let $s\in G$ be arbitrary, and
consider the functional $1_H.\lambda(s)$ 
defined by $1_H.\lambda(s)(a) = 1_H(\lambda(s)a)$, $a\in\cstar(G)$.
Since $F_\perp$ is an ideal also in $\M(\cstar(G))$, 
it follows that $1_H.\lambda(s) \in F$ because $F$ is weak* closed. 
Actually $1_H.\lambda(s)$ is just the left translation 
$\lt_s 1_H$ so by Lemma~\ref{lemma:const-on-cosets}
\[
1_H(shs\inv) = \lt_s 1_H(hs\inv) = \lt_s 1_H(s\inv) = 1
\]
for every $h$ in $H$. We conclude that $H$ is normal.

To show the final statement, suppose that $H$ is normal
and write $(\H,\pi)$ for the compact quantum subgroup
induced by $X = \cstar(H)$.  
The map $Q\col \lone(G)\to \lone(G/H)$ defined by
\[
Q f(sH) = \int_H f(sh) \, dh \qquad(f\in\lone(G),\, s\in G)
\] 
is a $*$-homomorphism, so it induces a $*$-homomorphism 
$\rho\col \cstar(G)\to \cstar(G/H)$. 
The pair $((G/H)\latedual,\rho)$
is a compact quantum subgroup of $\G = \dual G$. 
A simple calculation shows that the 
adjoint map $\rho^*\col \B(G/H)\to \B(G)$ satisfies 
$\rho^*u(s) = u(sH)$ for every $u$ in $B(G/H)$ and $s$ in $G$. 
Therefore $\rho^*$ is an isometry from $\B(G/H)$ onto $F$
by \cite[Corollaire~2.26]{eymard:fourier}.
It follows that $F_\perp = \ker \rho$. We conclude that 
$(\H, \pi)\cong((G/H)\latedual,\rho)$.
\end{proof}

The preceding theorem is connected to the following result due to Bekka, Lau 
and Schlichting \cite[Corollary~1.4]{bekka-lau-schlichting}:
there is a one-to-one correspondence between closed, normal subgroups 
of an amenable locally compact group $G$ and weak*-closed, 
translation invariant $*$-subalgebras of the Fourier--Stieltjes algebra 
$\B(G)$.

The next result follows from Theorems~\ref{thm:subgroup->subalgebra}, 
\ref{thm:dual-of-LL}, \ref{thm:dual-normal} and~\ref{thm:2nd-unique}.

\begin{corollary}
Suppose that $G$ is an amenable locally compact group.
If $\H$ is a compact quantum subgroup of\/ $\G = \dual G$, then there
is an open, normal subgroup $H$ of $G$ such that 
$(\H, \pi)\cong ((G/H)\latedual,\rho)$.  
\end{corollary}

\begin{acknow}
I thank the referee for several insightful comments and 
suggestions, which have greatly improved the paper. 
In particular, the approach using $X$-trivial
representations is due to the referee, who also 
suggested the much shorter proof of the implication
(\ref{item:H-normal})$\implies$(\ref{item:symmetric})
in Theorem~\ref{thm:dual-normal}.
The $X$-trivial approach extends the original construction
to the non-symmetric case and also 
by-passes a gap noticed by the referee
(the gap could have been filled otherwise too).
Both Adam Skalski and the referee deserve thanks for  
pointing out the very relevant reference \cite{tomatsu:coideal};
the referee also for references 
\cite{vaes:imprimitivity,enock:sous-facteur}.
I also thank Zhiguo Hu, Adam Skalski and Nico Spronk for comments. 
\end{acknow}

\end{document}